\documentclass[10pt,a4paper]{amsart}
\usepackage[foot]{amsaddr}
\usepackage{amsmath}
\usepackage{amsthm}
\usepackage{amssymb}
\usepackage{mathrsfs}
\usepackage{verbatim}
\usepackage{xcolor}
\usepackage[colorlinks=true,linkcolor=blue,citecolor=red]{hyperref}
\usepackage{enumitem}
\usepackage{multicol}
\usepackage{graphicx}
\usepackage{float}
\usepackage[nameinlink]{cleveref}

\theoremstyle{plain}
\newtheorem{theorem}{Theorem}
\newtheorem*{theorem*}{Theorem}
\newtheorem{lemma}[theorem]{Lemma}
\newtheorem{proposition}[theorem]{Proposition}

\theoremstyle{definition}

\theoremstyle{remark}
\newtheorem{remark}[theorem]{Remark}
\newtheorem{example}[theorem]{Example}

\newcommand{\RR}{\mathbb{R}}
\newcommand{\CC}{\mathbb{C}}
\newcommand{\NN}{\mathbb{N}}

\newcommand{\test}{C_c^{\infty}}

\newcommand{\Om}{\Omega}
\newcommand{\II}{\iint\limits_}
\newcommand{\I}{\int\limits_}

\numberwithin{equation}{section}
\numberwithin{theorem}{section}
\allowdisplaybreaks

\makeatletter
\@namedef{subjclassname@2020}{\textup{2020} Mathematics Subject Classification}
\makeatother

\begin{document}
	
	\title[Magnetic fractional Poincar\'e inequality in punctured domains]{Magnetic fractional Poincar\'e inequality in punctured domains}
	
	\author{$^1$Kaushik Bal}
	
	\email{$^1$kaushik@iitk.ac.in}
	
	\author{$^2$Kaushik Mohanta}
	\email{$^2$kaushik.k.mohanta@jyu.fi}
	\author{$^3$Prosenjit Roy}
	
	\email{$^3$prosenjit@iitk.ac.in}
	
	\address{$^{1,3}$Department of Mathematics and Statistics, Indian Institute of Technology Kanpur, India}
	\address{$^{2}$Department of Mathematics and Statistics, University of Jyv\"askyl\"a, Finland}

	.
	\smallskip
	\begin{abstract}
		We study Poincar\'e-Wirtinger type inequalities in the framework of magnetic fractional Sobolev spaces. In the local case, Lieb-Seiringer-Yngvason [E. Lieb, R. Seiringer, and J. Yngvason, Poincar\'e inequalities in punctured domains, Ann. of Math., 2003] showed that, if a bounded domain $\Om$ is the union of two disjoint sets $\Gamma$ and $\Lambda$, then the $L^p$-norm of a function calculated on $\Om$ is dominated by the sum of magnetic seminorms of the function, calculated on $\Gamma$ and $\Lambda$ separately. We show that the straightforward generalisation of their result to nonlocal setup does not hold true in general. We provide an alternative formulation of the problem for the nonlocal case. As an auxiliary result, we also show that the set of eigenvalues of the magnetic fractional Laplacian is discrete.
	\end{abstract}
	
	\maketitle


	\section{Introduction}
	The objective of this article is to study Poincar\'e type  inequalities in punctured domains for magnetic fractional Sobolev spaces. By Poincar\'e inequality in a domain $\Om\subseteq \RR^N$, we generally understand an inequality of the form
	\begin{equation}\label{first-poincare}
	\| f\|_{L^2(\Om,\RR)}\leq C\| \nabla f\|_{L^2(\Om,\RR)},
	\end{equation}
	where $C=C(N,\Om)>0$ is a constant. Clearly \eqref{first-poincare} fails when $f:\Om\to\RR$ is a non-zero constant function and $\Om$ is bounded, but when we restrict $f$ to certain subclasses of $W^{1,2}(\Om,\RR)$ in order to ``keep them away" from the one dimensional subspace of constant functions, the inequality does hold. There are two commonly used subspaces which are considered for this purpose:  $\test(\Om,\RR)$, related to the Dirichlet eigenvalue problems (see \cite[Chapter~5.8.1]{evans}) and the subspace of functions with with the property $\frac{1}{\mathcal{L}^N(\Om)}\I{\Om}f=0$, related to the Neumann eigenvalue problem (see \cite[Theorem~8.11]{LiLo}
	. In the later case, the inequality is also referred to as Poincar\'e-Wirtinger inequality. This can be generalised further by taking $p\geq 1$ and $1\leq q\leq \frac{Np}{N-p}$ to be the exponents of the right and left hand sides respectively in \eqref{first-poincare}. It is then called Poincar\'e-Wirtinger-Sobolev inequality. Regarding this, we have the following theorem, known in literature: 
	\begin{theorem}[\protect{\cite[Theorem~8.11 and Theorem~8.12]{LiLo}}] \label{lieb-th1}
		Let $\Om\ (\subseteq\RR^N)$ be a bounded domain with cone property. Let $q\in[1,\infty]$ and $p\in\left[\max\left\{ 1,\frac{Nq}{N+q} \right\},\infty\right]$ if $q<\infty$, $p\in(N,\infty]$ if $q=\infty$. Let $g\in L^{p'}(\Om,\CC)$, where $\frac{1}{p}+\frac{1}{p'}=1$, be such that $\I{\Om}g=1$. Then there exists a constant $S_{p,q}=S_{p,q}(\Om,g,p,q)>0$ such that for any $f\in W^{1,p}(\Om,\CC)$,
		\begin{equation*}
		\left\| f-\I{\Om}fg\right\|_{L^q(\Om,\CC)}\leq S_{p,q} \|\nabla f\|_{L^p(\Om,\CC)}.
		\end{equation*}
	\end{theorem}
	\bigskip
	
	In the theory of electromagnetism, the vector potential $A$ appears naturally as the vector field, up to translation by constants, for which $\mbox{curl}(A)=B$, where $B$ denotes the electromagnetic field. The operator $\nabla+iA$, defined by
	$$
	(\nabla + i A)f(x):=\nabla f(x)+iA(x)f(x),
	$$
	plays an important role in quantum mechanics, specifically in the study of Bose-Einstein condensations (BEC) and superfluidity. The map $f\mapsto \|(\nabla + i A) f\|_{L^p(\Om,\CC)}$, where $A$ is a vector field on $\Om$, gives a seminorm on the Sobolev spaces, which can be regarded as generalisations of the usual Sobolev seminorms: $f\mapsto \|\nabla f\|_{L^p(\Om,\CC)}$. Various generalisations of the standard Poincar\'e inequality serve as major tools in proving results in BEC and superfluidity, \cite{LiSe,LiSeYn02}; some require the splitting of the domains while others require replacing the role of the operator $\Delta$ by $\Delta+i A$, $A$ being the so called magnetic potential. For a detailed discussion on this, the reader is referred to \cite{LiSeSoYn}.
	
	In \cite{LiSeYn}, some remarkable results in this direction were proved by Lieb-Seiringer-Yngvason. They generalised \Cref{lieb-th1} in two ways: 
	The first generalisation comes through replacing the operator $\nabla$ by $\nabla + i A$. It is given as follows:
	
	\begin{theorem}[\protect{\cite[Theorem~2]{LiSeYn}}]\label{lieb-th2}
		Let $\Om\ (\subseteq\RR^N)$ be a bounded domain with cone property and $A:\Om\to\RR^N$ be a bounded vector field on $ \Om$. Let $q\in[1,\infty]$. Set $r=\max\left\{1,\frac{qN}{N+q}\right\}$ if $q<\infty$ and $r>N$ if $q=\infty$. Take $p\in(r,\infty]$. Define the energy
		\begin{equation*}
		E^{p,q}_{A}:=\inf \left\{ \frac{\|(\nabla +iA)f(x)\|_{L^p(\Om,\CC)}}{\| f \|_{L^q(\Om,\CC)}} \ \Big| \ f\in W^{1,p}(\Om,\CC)\setminus \{0\}\right\},
		\end{equation*}
		and the ground state manifold
		\begin{equation*}
		M^{p,q}_{A}:= \left\{ f\in W^{1,p}(\Om,\CC) \ \Big| \ \frac{\|(\nabla +iA)f(x)\|_{L^p(\Om,\CC)}}{\| f \|_{L^q(\Om,\CC)}}=E^{p,q}_A \right\}.
		\end{equation*}
		Let $0<\delta\leq 1$. 
		Then there exists a constant $S^{p,q}_\delta>0$ such that for any $f\in W^{1,p}(\Om,\CC)$ with $\inf\limits_{\phi\in M^{p,q}_A}\|f-\phi\|_{L^q(\Om,\CC)}\geq \delta \|f\|_{L^q(\Om,\CC)}$, one has
		\begin{equation*}
		\|(\nabla+iA)f\|_{L^p(\Om,\CC)}\geq \left(\frac{1}{S^{p,q}_\delta}+E^{p,q}_A\right)\|f\|_{L^q(\Om,\CC)}.
		\end{equation*}
	\end{theorem}
	\bigskip
	
	Observe that \Cref{lieb-th2} is indeed a generalisation of \Cref{lieb-th1}, at least in the particular case, when $g\equiv \frac{1}{\mathcal{L}^N(\Om)}$.
\begin{remark}\label{lieb implies usual}
	Let $p,q,\Om$ as in \Cref{lieb-th2} with $p,q< \infty$, $A\equiv0$, for an arbitrary $f_1\in L^q(\Om,\CC)$, denote $f:=f_1-\frac{1}{\mathcal{L}^N(\Om)}\I{\Om}f_1$. We have, for any $c\in\CC$,
	\begin{equation*}
		|c\mathcal{L}^N(\Om)|= |\I{\Om}(c-f)|\leq\I{\Om}|f-c|\leq \|f-c\|_{L^q(\Om,\CC)}(\mathcal{L}^N(\Om))^\frac{1}{q'},
	\end{equation*}
	where we used Holder's inequality. This gives
	\begin{equation*}
		|c|\mathcal{L}^N(\Om)^\frac{1}{q}\leq  \|f-c\|_{L^q(\Om,\CC)}.
	\end{equation*}
	We use this inequality to obtain
	\begin{align*}
		\|f\|_{L^q(\Om,\CC)}
		&\leq \|f-c\|_{L^q(\Om,\CC)}+\|c\|_{L^q(\Om,\CC)}\\
		&= \|f-c\|_{L^q(\Om,\CC)}+|c|\mathcal{L}^N(\Om)^\frac{1}{q}\\ 
		&\leq  2 \|f-c\|_{L^q(\Om,\CC)}.
	\end{align*}
	This shows that $\inf\limits_{\phi\in M^{p,q}_A}\|f-\phi\|_{L^q(\Om,\CC)}\geq \frac{1}{2} \|f\|_{L^q(\Om,\CC)}$, as in the case $A\equiv 0$, $M^{p,q}_A$ contains only constant functions. If we assume \Cref{lieb-th2} to be true, we have
	\begin{equation*}
		\left\| f_1-\frac{1}{\mathcal{L}^N(\Om)}\I{\Om}f_1\right\|_{L^q(\Om,\CC)}\leq C \|\nabla f_1\|_{L^p(\Om,\CC)}.
	\end{equation*}
	So, \Cref{lieb-th2} is stronger than \Cref{lieb-th1}.
\end{remark}
\smallskip 
		
	The second generalisation (as mentioned before \Cref{lieb-th2}) deals with the case of ``punctured domains''.
	
	\begin{theorem}[\protect{\cite[Theorem~3]{LiSeYn}}]\label{lieb-th3}
		Let $\Om$, $p$, $q$, $r$, $E^{p,q}_{A}$ and $M^{p,q}_{A}$ be as in \Cref{lieb-th2}.
		Let $\Lambda\subseteq\Om$ be measurable,  $\Gamma:=\Om\setminus\Lambda$ and take $0<\delta\leq 1$. 
		Then for any $\varepsilon>0$, there exists $C=C(\Om,A,p,q,\delta,\varepsilon)>0$ such that for any $f\in W^{1,p}(\Om,\CC)$ with $\inf\limits_{\phi\in M^{p,q}_A}\|f-\phi\|_{L^q(\Om,\CC)}\geq \delta \|f\|_{L^q(\Om,\CC)}$,
		\begin{equation*}
		\|(\nabla+iA)f\|_{L^p(\Lambda,\CC)}+C\|(\nabla+iA)f\|_{L^r(\Gamma,\CC)}\geq \left(\frac{1}{S^{p,q}_\delta+\varepsilon}+E^{p,q}_A\right)\|f\|_{L^q(\Om,\CC)}
		\end{equation*}
		where $S^{p,q}_\delta$ is the optimal constant in \Cref{lieb-th2}. 
	\end{theorem}
	
	
	Our aim in this paper is to present the appropriate nonlocal analogues of \Cref{lieb-th2,lieb-th3}. Before proceeding further, let us first introduce our functional setup, which is the foundation for the work done in this article. For fixed open subset $\Om\subseteq\RR^N$, $s\in(0,1)$, a vector field $A:\RR^N\to\RR^N$, and $p\geq 1$, we define the \emph{magnetic fractional Sobolev space}
	\begin{equation*}
	W^{s,p}_A(\Om,\CC):=\{f\in L^p(\Om,\CC)\ \Big| \ [f]_{W^{s,p}_A(\Om,\CC)}<\infty\},
	\end{equation*}
	where \begin{equation*}
	[f]_{W^{s,p}_A(\Om,\CC)}:=\left(\ \II{\Om\times\Om}\frac{|f(x)-e^{i(x-y)\cdot A\left(\frac{x+y}{2}\right)}f(y)|^p}{|x-y|^{N+sp}}dxdy\right)^\frac{1}{p}
	\end{equation*}
	gives a seminorm on $W^{s,p}_A(\Om,\CC)$, the norm on the space is defined to be
	\begin{equation*}
	\| f \|_{W^{s,p}_A(\Om,\CC)}:=\left(\| f \|_{L^p(\Om,\CC)}^p+[ f ]_{W^{s,p}_A(\Om,\CC)}^p\right)^\frac{1}{p}.
	\end{equation*}
	When $A\equiv 0$, we denote the above space, seminorm and norm by $W^{s,p}(\Om,\CC)$, $[\cdot]_{{W^{s,p}(\Om,\CC)}}$ and $\|\cdot \|_{{W^{s,p}(\Om,\CC)}}$ respectively. In \cite{SqVo}, it has been shown that the Bourgain-Brezis-Mironescu formula holds for these spaces; that is with appropriate scaling the magnetic fractional seminorm converges to the magnetic seminorm. Thus these spaces can be viewed as the nonlocal analogues of the magnetic Sobolev spaces. Some fundamental results regarding fractional magnetic Sobolev spaces were studied in \cite{dASq,PiSqVe,BoSa21,MaSaVe}. The special case where $A\equiv0$ is the well-known \textit{fractional Sobolev space} ${W^{s,p}(\Om,\CC)}$. For general results regarding ${W^{s,p}(\Om,\RR)}$, we refer the reader to \cite{hhg} and the references therein.
	
	To the best of our knowledge, fractional Poincar\'e-Wirtinger inequality was first established by Ponce \cite{Po} as an application of BBM formula (see \cite{BBM}). In \cite{HuVa} a more general version of the result was proved. In the fractional magnetic Sobolev setup some Sobolev and Hardy-Sobolev type inequalities were studied in \cite{GuMe,LiChGu,LiJiGu}. We shall prove the following result which, we believe, is known to the experts; the local variant of this result is well known (see \cite[Chapter~8.8]{LiLo}).
	
	\begin{theorem}[Fractional Poincar\'e-Wirtinger-Sobolev inequality]\label{std-poincare}
		Let $\Om\subseteq\RR^N$ be a bounded domain with Lipschitz boundary, $s\in(0,1)$, $p\in [1,\infty)$, $q \in [1,\infty]$, $g\in L^{p'}(\Om)$ with $\I{\Om}g=1$, where $\frac{1}{p}+\frac{1}{p'}=1$. Further, assume that one of the following three conditions hold:\\
		(i) $sp<N$ and $q\leq\frac{Np}{N-sp}$,\\
		(ii) $sp=N$ and $q<\infty$,\\
		(iii) $sp>N$ and $q\leq\infty$.\\
		
		\noindent Then there exists a constant $C=C(\Om,g,p,q,s)>0$ such that for any $f\in W^{s,p}(\Om,\CC)$
		\begin{equation*}
		\left\| f-\I{\Om}fg \right\|_{L^q(\Om,\CC)}\leq C[f]_{W^{s,p}(\Om,\CC)}.
		\end{equation*}
	\end{theorem}
	\smallskip

	Let us now state our main results, after introducing some terminologies required for formulating the statements. For $s\in(0,1)$, $p\in [1,\infty)$ and $q\in [1,\infty]$ we define the energy
	\begin{equation*}
	E^{p,q}_{s,A}:=\inf \left\{ \frac{[f]_{W^{s,p}_A(\Om,\CC)}}{\| f \|_{L^q(\Om,\CC)}} \ \Big| \ f\in W^{s,p}(\Om,\CC),f\neq 0 \right\}.
	\end{equation*}
	The corresponding ground state manifold  is defined to be
	\begin{equation*}
	M^{p,q}_{s,A}:=\left\{ f\in W^{s,p}(\Om,\CC)\ \Big|\  \frac{[f]_{W^{s,p}_A(\Om,\CC)}}{\| f \|_{L^q(\Om,\CC)}}=E^{p,q}_{s,A}\right\}.
	\end{equation*}
	We use the following notion of distance from the ground state manifold:
	\begin{equation*}
	d^q_{s,A}(f):=\inf\limits_{\phi\in M^{p,q}_{s,A}}\| f-\phi\|_{L^q(\Om,\CC)}.
	\end{equation*}

	\noindent Our first main result is a generalisation of \Cref{lieb-th2} to the nonlocal case. 
	\begin{theorem}\label{th-poincare-for-vf}
		Let $\Om\ (\subseteq\RR^N)$ be a bounded Lipschitz domain, $A$ be a bounded vector field on the convex hull of $\Om$ and $s\in(0,1)$. Assume that $p\in[1,\infty)$ and $q\in[1,\infty]$ satisfy one of the following three conditions:\\
		(i) $sp<N$ and $q<\frac{Np}{N-sp}$,\\
		(ii) $sp=N$ and $q<\infty$,\\
		(iii) $sp>N$ and $q\leq\infty$.\\
		Fix $\delta\in(0,1] $. Then there is a constant $S=S(p,q,s,\delta,\Om)>0$ such that $\forall f\in W^{s,p}(\Om,\CC)$ with $d^q_{s,A}(f)\geq \delta\| f \|_{L^q(\Om,\CC)}$,
		\begin{equation}\label{th-poincare-for-vf-eq1}
		\| f \|_{L^q(\Om,\CC)}\leq S \left([f]_{W^{s,p}_A(\Om,\CC)}-E^{p,q}_{s,A}\|f\|_{L^q(\Om,\CC)}\right).
		\end{equation}
		The best constant $S$ in the above inequality is achieved. Consequently the ground state manifold $M^{p,q}_{s,A}$ is nonempty.
	\end{theorem}
	\smallskip
	
	Our second main result is a generalisation of \Cref{lieb-th3} in the nonlocal setup. Before stating it let us make a small but important observation which confirms that one cannot expect a straightforward generalisation of \Cref{lieb-th3}. To see this, we provide the following example:
	\begin{example}\label{example-1}
		We shall deal with the case $A\equiv 0$ here. For $p,q\geq 1,\ s\in(0,\frac{1}{p})$, take $\Om=B(0,1),\ \Lambda=B(0,\frac{1}{2}),\ \Gamma=B(0,1)\setminus B(0,\frac{1}{2})$. Take $f=\chi_{\Lambda}\in W^{s,p}(\Om,\CC)$ (see \cref{ap-char}). Then $[f]_{W^{s,p}(\Lambda,\CC)}=[f]_{W^{s,p}(\Gamma,\CC)}=0$ but $\|f\|_{L^q(\Om,\CC)}>0$.
		So, no matter what is the value of $C$,
		$$
		\frac{[f]_{W^{s,p}(\Lambda,\CC)}+C[f]_{W^{s,r}(\Gamma,\CC)}}{\|f\|_{L^q(\Om,\CC)}}=0.
		$$
		Hence an inequality of the form 
		\begin{equation}\label{example-inequality}
		[f]_{W^{s,p}(\Lambda,\CC)}+C[f]_{W^{s,r}(\Gamma,\CC)}\geq \left(\frac{1}{S+\varepsilon}+E^{p,q}_{s,A}\right)\|f\|_{L^q(\Om,\CC)}
		\end{equation}
		can not be expected to hold in the nonlocal setup.
	\end{example}
	\noindent However one may ask what happens if we impose an additional condition: $sp>1$, which is the required condition for the Dirichlet-type Poincar\'e inequality to hold in bounded domains (see \cite{hhg,Ind,MoSk,BMRS} for details). Here the scenario does not change, as we can see from the following example. Before going to the example note that in the hypotheses of \Cref{lieb-th3}, we assumed $1\leq r<p$.
	\begin{example}\label{example-2}
		As before we are in the case $A\equiv 0$. Note that in this case, $E^{p,q}_{s,A}=0$ and $M^{p,q}_{s,A}$ is the one-dimensional subspace consisting of constant functions. Let $\Om=(-1,1)\times(0,1)(\subseteq \RR^2)$, $q\in (1,\infty), r\in[1,p)$. Fix $s\in(\frac{1}{p},\frac{1}{r})$. Set $\Lambda:=(-1,0]\times (0,1)$ and $\Gamma:=(0,1)\times (0,1)$. For $\varepsilon>0$, we define $f_\varepsilon: \Om\to\RR$ by the formula 
		
		\begin{equation}\label{example-function}
		f_\varepsilon (x_1,x_2):=\begin{cases}
		\frac{2-3\varepsilon}{2+\varepsilon},\quad -1<x\leq0\\
		\frac{2-3\varepsilon}{2+\varepsilon}-\frac{2(2-\varepsilon)}{\varepsilon(2+\varepsilon)}x_1,\quad 0<x\leq\varepsilon\\
		-1,\quad \varepsilon<x<1.
		\end{cases}
		\end{equation}
		Note that $\I{\Om}f_\varepsilon(x)dx=0$, $ f_\varepsilon \in W^{1,p}(\Om,\CC)\subseteq W^{s,p}(\Om,\CC)$,  $[f_\varepsilon]_{W^{s,p}(\Lambda,\CC)}=0$ for any $\varepsilon>0$, and $\|f_\varepsilon \|_{L^q(\Om,\CC)}\to 1$ as $\varepsilon\to 0$. Moreover $[f_\varepsilon]_{W^{s,r}(\Gamma,\CC)}\to 0$ as $\varepsilon\to 0$ (see \Cref{appendix1} for details). Hence the left hand side of the nonlocal analogue of the inequality \eqref{example-inequality} cannot just consist of the seminorms of the two components, even when we have $sp>1$.
	\end{example}

	\noindent Further, \Cref{example-1,example-2} give us a hint on what should be the appropriate fractional analogue of \Cref{lieb-th3}, as it appears that, in the nonlocal case, the terms coming from integrating over $\Gamma\times\Lambda$ is non-negligible, and we must take them into account. This is clarified in our second main result:
	
	\begin{theorem}\label{th-poincare-punctured}
		Let $\Om\ (\subseteq \RR^N)$ be a bounded open set with Lipschitz boundary, $A$ be a bounded vector field on the convex hull of $\Om$, $s\in(0,1)$ and $\Lambda\subseteq\Om$. Assume $1\leq q\leq\infty$ and $1\leq r<p<\infty$ are such that one of the following conditions holds:\\
		(a) $r=1$ and $q<\frac{N}{N-s}$,\\
		(b) $sr>N$ and $q=\infty$,\\
		(c) $q<\frac{Nr}{N-sr}$, $sr<N$ and $\frac{N}{N-s}\leq q<\infty$;\\
		and that one of the following conditions hold:\\
		(i) $sp<N$ and $q<\frac{Np}{N-sp}$,\\
		(ii) $sp=N$ and $q<\infty$,\\
		(iii) $sp>N$ and $q\leq\infty$.\\		
		Fix $\delta\in(0,1]$. For any $\varepsilon>0$ there is some $C=C(\Om,A,s,p,q,\delta,\varepsilon)>0$ such that for all $f\in W^{s,p}(\Om,\CC)$ with $d^q_{s,A}(f)\geq\delta\|f\|_{L^q(\Om,\CC)}$,
		
		\begin{align*}
		&\left(\ \II{\Lambda\times\Lambda}\frac{|f(x)-e^{i(x-y)\cdot A\left(\frac{x+y}{2}\right)}f(y)|^p}{|x-y|^{N+sp}}dxdy\ \right )^\frac{1}{p}\\
		&+C\left(\ \II{(\Om\times\Om)\setminus(\Lambda\times\Lambda)}\frac{|f(x)-e^{i(x-y)\cdot A\left(\frac{x+y}{2}\right)}f(y)|^r}{|x-y|^{N+sr}}dxdy\ \right)^\frac{1}{r}\\
		&\geq\left(\frac{1}{S+\varepsilon}+E^{p,q}_{s,A}\right)\|f\|_{L^q(\Om,\CC)}
		\end{align*}
		where $S$ is the best constant in \Cref{th-poincare-for-vf}.
	\end{theorem}

	Our technique of proving the \Cref{th-poincare-for-vf,th-poincare-punctured} is motivated by the techniques used in \cite{LiSeYn}, but is significantly different due to many reasons, one of which is the non-availability of embeddings of fractional Sobolev spaces of different exponents. Mironescu-Sickel \cite{MiSi} have shown that, unlike in the case $s=1$, when $s\in(0,1)$, $1\leq r<p<\infty$ and $\Om$ is a bounded domain, the embedding $W^{s,p}(\Om,\RR)\subseteq W^{s,r}(\Om,\RR)$ never hold. However, we have shown in \Cref{p-r-embedding} that if we take $s_1<s_2$, $W^{s_2,p}(\Om,\CC)\subseteq W^{s_1,r}(\Om,\CC)$. Similar results can also be found in \cite{AnWa}. However we provide an independent proof of this important result. \Cref{p-r-embedding} is used to overcome the non-availability of the embedding $W^{s_2,p}(\Om,\CC)\subseteq W^{s_1,r}(\Om,\CC)$ in the Step-2/b in the proof of \Cref{th-poincare-punctured}. However, in the process we loose information about the case $sr<N$ with $q=\frac{Nq}{N-sq}$.\smallskip
	
	The magnetic fractional Sobolev spaces introduced above are the ideal spaces to study problems related to the so called regional magnetic fractional $p$-Laplacian $(-\Delta_{p,A})^s$, which is defined by the formula:
	\begin{equation*}
		(-\Delta_{p,A})^sf(x)
		:=pv \I{\Om} \frac{|f(x)-e^{i(x-y)\cdot A\left(\frac{x+y}{2}\right)}f(y)|^{p-2}(f(x)-e^{i(x-y)\cdot A\left(\frac{x+y}{2}\right)}f(y))}{|x-y|^{N+sp}}dy.
	\end{equation*}
	In the special case $p=2$, it is called the \textit{regional magnetic fractional Laplacian}, and it is the operator, on $W^{s,2}(\Om,\CC)$, associated with the quadratic form $[\cdot]_{W^{s,2}_A(\Om,\CC)}$. For the classical case $A\equiv 0$, this operator has been studied in recent works, for example, see \cite{AbFaTe,Fall,Guan,Temgoua,Warma16}. The operator $(-\Delta_{2,A})^s$ can be compared to the \textit{magnetic fractional Laplacian}, whose definition is similar to that of $(-\Delta_{2,A})^s$ as above, but the domain of integration being $\RR^N$. Some recent works on this operator, for non-trivial $A$, can be found in \cite{Ambrosio22,AmdAv,ZuLo}. In \Cref{laplacian}, we study the operator $(-\Delta_{2,A})^s$ and prove that this is a self-adjoint operator with discrete spectrum. For $s=1$, similar results can be found in \cite{LiSeYn}.
	
	The article is organised in the following way: 
	in \Cref{preli} we recall some known results, which we shall use in the later sections. In \Cref{main} we prove \Cref{std-poincare,,th-poincare-for-vf,,th-poincare-punctured}. In \Cref{laplacian}, we introduce the magnetic fractional Laplacian and discuss some of its basic properties. In \Cref{discrete}, we prove \Cref{self adjoint,discrete} regarding the regional magnetic fractional Laplacian.
	\smallskip

	\section{Some preliminary results}\label{preli}
	In this section we recall some results, which, up to some small modifications, are already known in literature. Throughout this article we shall use the following notations:
	\begin{enumerate}
		\item $\mathcal{L}^N$ will denote the $N$-dimensional Lebesgue/Hausdorff measure.
		\item $C$ will stand for an arbitrary constant, which can change from line to line.
		\item $\mathbb{S}^{N-1}$ stands for the $N-1$-dimensional sphere, centred at the origin.
		\item The notation $f_n\rightharpoonup f$ would mean that $f_n$ converges to $f$ weakly (in some space to be specified).
	\end{enumerate}

	\begin{lemma}[An embedding result]\label{p-r-embedding}
		Let $0<s_1<s_2<1$, $1\leq r\leq p<\infty$, $\Om\subseteq \RR^N$ and $H\subseteq \Om\times\Om$. Let $A$ be a bounded vector field on the convex hull of $\Om$. Then for any $f\in L^p(\Om,\CC)$ 
		\begin{align*}
		\II{H}&\frac{|f(x)-e^{i(x-y)\cdot A\left(\frac{x+y}{2}\right)}f(y)|^r}{|x-y|^{N+s_1r}}dxdy\\	
		\leq& \left(\II{H}\frac{|f(x)-e^{i(x-y)\cdot A\left(\frac{x+y}{2}\right)}f(y)|^p}{|x-y|^{N+s_2p}}dxdy\right)^\frac{r}{p}
		\left(\II{H}|x-y|^\frac{Nr+s_2rp-Np-s_1rp}{p-r}dxdy \right)^\frac{p-r}{p}.
		\end{align*}
		Moreover if $\Om$ is bounded, then for some $C=C(\Om,N,p,r)>0$, 
		$$
		[f]_{W^{s_1,r}_A(\Om,\CC)}^r\leq \frac{C(\Om,N,p,r)}{(s_2-s_1)^\frac{p-r}{r}}[f]_{W^{s_2,p}_A(\Om,\CC)}^r.
		$$
	\end{lemma}
	\begin{proof}
		The case $p=r$ is similar to the remaining part of the proof, also it is well known (see \cite{hhg}), hence the proof is omitted. So we assume that $\frac{p}{r}>1$ and observe that $\frac{1}{\frac{p}{r}}+\frac{1}{\frac{p}{p-r}}=1$.
		\begin{align*}
		\II{H}&\frac{|f(x)-e^{i(x-y)\cdot A\left(\frac{x+y}{2}\right)}f(y)|^r}{|x-y|^{N+s_1r}}dxdy\\
		=&\II{H}\frac{|f(x)-e^{i(x-y)\cdot A\left(\frac{x+y}{2}\right)}f(y)|^r}{|x-y|^{\frac{Nr}{p}+s_2r}}|x-y|^{\frac{Nr}{p}+s_2r-N-s_1r}dxdy\\
		\leq& \left(\II{H}\frac{|f(x)-e^{i(x-y)\cdot A\left(\frac{x+y}{2}\right)}f(y)|^p}{|x-y|^{N+s_2p}}dxdy\right)^\frac{r}{p}
		\left(\II{H}|x-y|^\frac{Nr+s_2rp-Np-s_1rp}{p-r}dxdy \right)^\frac{p-r}{p}.
		\end{align*}
		Therefore, the first part is proved. Now for the second part, assume $\mbox{diam}(\Om)<R$ and put $H=\Om\times\Om$, in the above calculation, to get
		\begin{align*}
		[f]_{W^{s_1,r}_A(\Om,\CC)}^r
		\leq& [f]_{W^{s_2,p}_A(\Om,\CC)}^r\left(\I{\Om}dy\I{B(0,R)}|t|^\frac{Nr+s_2rp-Np-s_1rp}{p-r}dt \right)^\frac{p-r}{p}\\
		=& \mathcal{L}^{N-1}(\mathbb{S}^{N-1})^\frac{p-r}{p} [f]_{W^{s_2,p}_A(\Om,\CC)}^r\mathcal{L}^N(\Om)^\frac{p-r}{p}
		\left(\I{0}^R t^\frac{Nr+s_2rp-Np-s_1rp}{p-r}t^{N-1}dt \right)^\frac{p-r}{p}\\
		=&\mathcal{L}^{N-1}(\mathbb{S}^{N-1})^\frac{p-r}{p} \mathcal{L}^N(\Om)^\frac{p-r}{p}[f]_{W^{s_2,p}_A(\Om,\CC)}^r
		\left(\I{0}^R t^\frac{(s_2-s_1)rp-p+r}{p-r}dt \right)^\frac{p-r}{p}\\
		=&\mathcal{L}^{N-1}(\mathbb{S}^{N-1})^\frac{p-r}{p}\mathcal{L}^N(\Om)^\frac{p-r}{p} [f]_{W^{s_2,p}_A(\Om,\CC)}^r
		\left(\frac{p-r}{(s_2-s_1)rp}R^\frac{(s_2-s_1)rp}{p-r} \right)^\frac{p-r}{p}\\
		\leq& \frac{C(\Om,N,p,r)}{(s_2-s_1)^\frac{p-r}{p}}[f]_{W^{s_2,p}_A(\Om,\CC)}^r.
		\end{align*}
	\end{proof}

	\begin{lemma}[\protect{Ascoli-Arzel\`a theorem, \cite[Theorem~4.4.8]{Kum}}]\label{ascoli}
		Let $X$ be a compact metric space. Let $C(X, \CC)$ be given the sup norm metric. Then a set $B\subseteq C(X,\CC)$ is compact if and only if $B$ is bounded, closed and equicontinuous.
	\end{lemma}

	 We have the following result:
	
	\begin{lemma}[Compact embedding]\label{cpt-embedding}
		Suppose $\Om\ (\subseteq\RR^N)$ is a bounded domain with Lipschitz boundary, $s\in(0,1)$, $p\in [1,\infty)$ and $q \in[1,\infty]$. Assume that one of the following three conditions hold:\\
		(i) $sp<N$ and $q<\frac{Np}{N-sp}$,\\
		(ii) $sp=N$ and $q<\infty$,\\
		(iii) $sp>N$ and $q=\infty$.\\
		Then any bounded sequence in $W^{s,p}(\Om,\CC)$ has a convergent subsequence in $L^q(\Om,\CC)$.
	\end{lemma}
	\begin{proof}[Sketch of the proof]~\linebreak
	The proof of the fact that $W^{s,p}(\Om,\RR)$ is compactly embedded in $L^q(\Om,\RR)$, when hypotheses (i) or (ii) hold, goes exactly in the same way as the proof of \cite[Theorem~9.16]{Bre11}. The main ingredients of this proof are continuous embedding results and Ascoli-Arzel\`a theorem. The nonlocal counterparts of the continuous embedding results can be found in \cite[Theorem~6.10,~7.1, Corollary~7.2, Theorem~8.2]{hhg}.
	
	It remains to check, whether we can conclude the same for $\CC$-valued function spaces. For this, note that if $\{f_n\}$ is a bounded sequence in $W^{s,p}(\Om,\CC)$, clearly $\mbox{Re}(f)$ and $\mbox{Im}(f)$ are bounded sequences in $W^{s,p}(\Om,\RR)$ and hence have convergent subsequences in $L^p(\Om,\RR)$. This implies that $\{f_n\}$ has a convergent subsequence in $L^p(\Om,\CC)$. So we can replace $\RR$-valued function spaces with $\CC$-valued spaces.
	\end{proof}

	\begin{lemma}[\protect{\cite[Theorem~6.5]{hhg}}]\label{sobolev-inequality-1}
		Let $s\in (0,1)$ and $p\in[1,\infty)$ be such that $sp < n$. Then there exists a positive constant $C=C(n,p,s)$ such that, for any $f\in W^{s,p}(\RR^N,\RR)$, we have
		
		\begin{equation*}
		\|f\|_{L^{\frac{Np}{N-sp}}(\RR^N,\RR)}\leq C [f]_{W^{s,p}(\RR^N,\RR)}.
		\end{equation*}
		
	\end{lemma}
	
	\begin{lemma}[Sobolev-type inequality]\label{sobolev-inequality}
		Let $\Om$ be a bounded domain with Lipschitz boundary, $s\in (0,1)$, $p\in[1,\infty)$ be such that $sp < n$. Then there exists a positive constant $C=C(N,p,s,\Om,g)$ such that, for any $f\in W^{s,p}(\Om,\CC)$, we have
		
		\begin{equation*}
		\|f\|_{L^{\frac{Np}{N-sp}}(\Om,\CC)}\leq C \|f\|_{W^{s,p}(\Om,\CC)}.
		\end{equation*}
		
	\end{lemma}
	\begin{proof}
		Let $P:W^{s,p}(\Om,\RR)\to W^{s,p}(\RR^N,\RR)$ be the extension operator. We apply \Cref{sobolev-inequality-1}, on $\mbox{Re}(f),\ \mbox{Im}(f)$, to get
		\begin{align*}
		\|f\|_{L^{\frac{Np}{N-sp}}(\Om,\CC)}^{\frac{Np}{N-sp}}
		=&\I{\Om}|f|^{\frac{Np}{N-sp}}
		=\I{\Om}\left( \mbox{Re}(f)^2+\mbox{Im}(f)^2 \right)^{\frac{Np}{2(N-sp)}}\\
		\leq& C(N,p,s) \left[ \I{\RR^N}|P(\mbox{Re}(f))|^{\frac{Np}{N-sp}}+\I{\RR^N}|P(\mbox{Im}(f))|^{\frac{Np}{N-sp}} \right]\\
		\leq& C(p,\Om,N,s) \left[ \|P(\mbox{Re}(f))\|_{W^{s,p}(\Om,\RR)}^p+\|P(\mbox{Im}(f))\|_{W^{s,p}(\Om,\RR)}^p \right]^{\frac{Np}{p(N-sp)}}\\
		\leq& C(p,\Om,N,s) \left[ \|\mbox{Re}(f)\|_{W^{s,p}(\Om,\RR)}^p+\|\mbox{Im}(f)\|_{W^{s,p}(\Om,\RR)}^p \right]^{\frac{Np}{p(N-sp)}}\\
		\leq& C(p,\Om,N,s) \|f\|_{W^{s,p}(\Om,\CC)}^\frac{Np}{(N-sp)}.
		\end{align*}
		Consequently, the proof follows.
	\end{proof}
	
	In order to prove a Poincar\'e inequality for magnetic fractional Sobolev spaces for a class of functions which do not depend on $A$, we just need to prove it for the case $A=0$ because we have the following diamagnetic inequality.
	\begin{lemma}[Pointwise diamagnetic inequality, \protect{\cite[Remark~3.2]{dASq}}]\label{th-diamagnetic}
		For any $x,y\in\RR^N$, $A:\RR^N\to\RR^N$ and for any $f:U\ (\subseteq\RR^N)\to\CC$ which is finite almost everywhere, we have
		\begin{equation*}
		\left||f(x)|-|f(y)|\right|\leq|e^{-i(x-y)\cdot A\left(\frac{x+y}{2}\right)}f(x)-f(y)|
		\end{equation*}
	\end{lemma}
	
	The next result for the special case $\Om=\RR^N$ can be found in \cite[Lemma~3.6]{BoSa21}; the proof in this case is similar. However, we give a sketch of the proof. 
	\begin{lemma}\label{equality-magnetic-usual}
		If $\Om\ (\subseteq\RR^N)$ is an open set, $A$ is a bounded vector field on the convex hull of $\Om$, $s\in(0,1)$ and $1\leq p<\infty$ then $W^{s,p}(\Om,\CC)=W^{s,p}_A(\Om,\CC)$. Moreover for some $C>0$
		\begin{equation}
		[f]_{W^{s,p}(\Om,\CC)}^p \leq C\left([f]_{W^{s,p}_A(\Om,\CC)}+\|f\|_{L^p(\Om,\CC)}\right)
		\end{equation}
		and
		\begin{equation}
		[f]_{W^{s,p}_A(\Om,\CC)}^p \leq C\left([f]_{W^{s,p}(\Om,\CC)}+\|f\|_{L^p(\Om,\CC)}\right)
		\end{equation}
		and therefore the two norms $\|.\|_{W^{s,p}(\Om,\CC)}$ and $\|.\|_{W^{s,p}_A(\Om,\CC)}$ are equivalent.
	\end{lemma}
	\begin{proof}
		Let $f\in W^{s,p}_A(\Om,\CC)$. Then
		\begin{align*}
		[f]_{W^{s,p}(\Om,\CC)}^p
		&=\II{\Om\times\Om}\frac{|f(x)-f(y)|^p}{|x-y|^{N+sp}}dxdy\\
		&=\II{\Om\times\Om}\frac{|f(x)-e^{i(x-y)\cdot A\left(\frac{x+y}{2}\right)}f(y)+e^{i(x-y)\cdot A\left(\frac{x+y}{2}\right)}f(y)-f(y)|^p}{|x-y|^{N+sp}}dxdy.
		\end{align*}
		Now the numerator in the integrand is dominated by
		$$
		2^{p-1}(|f(x)-e^{i(x-y)\cdot A\left(\frac{x+y}{2}\right)}f(y)|^p+|e^{i(x-y)\cdot A\left(\frac{x+y}{2}\right)}-1|^p|f(y)|^p).
		$$
		Hence we can write
		\begin{align}\label{cpaa-eq1}
		\nonumber[f]_{W^{s,p}(\Om,\CC)}^p
		\leq& 2^{p-1}[f]_{W^{s,p}_A(\Om,\CC)}^p \\
		& + 2^{p-1} \I{y\in\Om}|f(y)|^p \left[\ \I{|x-y|\geq 1} \frac{|e^{i(x-y)\cdot A\left(\frac{x+y}{2}\right)}-1|^p}{|x-y|^{N+sp}}dx
		+\I{|x-y|< 1} \frac{|e^{i(x-y)\cdot A\left(\frac{x+y}{2}\right)}-1|^p}{|x-y|^{N+sp}}dx \right]dy.
		\end{align}
		Note that the function $\theta\mapsto e^{i\theta}$ has Lipschitz constant $1$, and its image has diameter $2$. Using this and the Cauchy-Schwartz inequality, it can be shown that
		\begin{equation*}
		|e^{i(x-y)\cdot A\left(\frac{x+y}{2}\right)}-1|\leq \begin{cases}
		2, \ &|x-y|\geq 1,\\
		\| A \|_\infty |x-y|, \ &|x-y|<1.
		\end{cases}
		\end{equation*}
		Using this in \eqref{cpaa-eq1} it follows that
		
		\begin{align*}
		[f]_{W^{s,p}(\Om,\CC)}^p
		\leq& 2^{p-1} [f]_{W^{s,p}_A(\Om,\CC)}^p\\
		&+ 2^{p-1}\I{y\in\Om}|f(y)|^p \left[\ \I{|x-y|\geq 1} \frac{2^p}{|x-y|^{N+sp}}dx
		+\I{|x-y|< 1} \|A\|_\infty|x-y|^{p-N-sp}dx \right]dy\\
		=&2^{p-1}\left( [f]_{W^{s,p}_A(\Om,\CC)}^p
		+ \left[ \frac{2^p\mathcal{L}^{N-1}(S^{N-1})}{sp}+ \frac{\|A\|_\infty \mathcal{L}^{N-1}(S^{N-1})}{p-sp} \right]\|f\|_{L^p(\Om,\CC)}\right).
		\end{align*}
		The other inequality can be proved similarly.
	\end{proof}
	\smallskip

	\section{Main results}\label{main}
	
	We start this section by proving \Cref{std-poincare}.
	\begin{proof}[\textbf{Proof of \Cref{std-poincare}}]
		First, we shall assume that one of the conditions (i), (ii) (iii) in the hypothesis of \Cref{cpt-embedding} holds, so that we can apply the compact embedding results. Without loss of generality, we can assume that $\mathcal{L}^N(\Om)=1$. Note, that if we prove this result for $q\geq p$, then the proof for the case $q<p$ follows immediately as $\Om$ is bounded and we have the embedding of $L^q$-spaces via H\"older inequality. So, we assume $q\geq p$ and suppose the statement of the theorem is false. Then we shall get a sequence of non-zero functions $\{f_n\}_{n\geq 1}$ in $W^{s,p}(\Om,\CC)$, where we can assume, without loss of generality, $\|f_n-\I{\Om}f_ng\|_{L^q(\Om,\CC)}=1$ for all $n$, such that $[f_n]_{W^{s,p}(\Om,\CC)}\to 0$ as $n\to\infty$. Set $\psi_n=f_n-\I{\Om}f_ng$. By H\"older's inequality, $\|\psi_n\|_{L^p(\Om,\CC)}\leq 1$ is uniformly bounded, and hence $\|\psi_n\|_{W^{s,p}(\Om,\CC)}$ is uniformly bounded. So, we apply \Cref{cpt-embedding} to conclude that, up to a subsequence, $\psi_n$ converges to some $\psi$ in $L^q(\Om,\CC)$ and hence in $L^1(\Om,\CC)$. Clearly, 
		\begin{equation}\label{std-poincare-eq1}
		\|\psi\|_{L^q(\Om,\CC)}=1.
		\end{equation} Also,
		\begin{equation}\label{std-poincare-eq2}
		\left| \I{\Om}\psi g\right|
		=\left| \I{\Om}\psi_ng-\I{\Om}\psi g \right|
		\leq \I{\Om}|\psi_n-\psi||g|
		\leq \|\psi_n-\psi\|_{L^q}\|g\|_{L^{q'}}\to 0.
		\end{equation}
		Now since $\psi_n$ converges to $\psi$ strongly in $L^q(\Om,\CC)$, up to a subsequence, we can assume that $\psi_n\to \psi$ a.e. pointwise. Using Fatou's lemma, we get
		\begin{equation*}
		0=\lim\limits_{n\to \infty}[\psi_n]_{W^{s,p}_A(\Om,\CC)}\geq [\psi]_{W^{s,p}_A(\Om,\CC)}.
		\end{equation*}
		This implies $\psi$ is a constant almost everywhere in $\Om$; furthermore, \eqref{std-poincare-eq1} implies that this constant is not $0$. This and \eqref{std-poincare-eq2} implies $\I{\Om}g=0$. This is a contradiction to the hypothesis on $g$.
        
        Next we consider the remaining case, that is $sp<N$ with $q=\frac{Np}{N-sp}$. The result follows from \Cref{sobolev-inequality} and the above case, as shown in the following calculation:
        \begin{align*}
        \left\| f-\I{\Om}fg \right\|_{L^q(\Om,\CC)}
        \leq& C\left\| f-\I{\Om}fg \right\|_{W^{s,p}(\Om,\CC)}\\
        =&C\left([f]_{W^{s,p}(\Om,\CC)}^p+\left\| f-\I{\Om}fg \right\|_{L^p(\Om,\CC)}\right)^\frac{1}{p}
        \leq C [f]_{W^{s,p}(\Om,\CC)}.        
        \end{align*}
	\end{proof}

	\begin{lemma}\label{seminorm-bdd-norm}
		Let $\Om$ be a bounded domain in $\RR^N$ with Lipschitz boundary, $A$ be a bounded vector field in the convex hull of $\Om$, $p,q\geq 1$, $s\in(0,1)$. Assume that for a sequence of functions $f_n$, $\{\|f_n\|_{L^q(\Om,\CC)}\}_n$ and $\{[f_n]_{W^{s,p}_A(\Om,\CC)}\}_n$ are bounded sequences. Then $\{\|f_n\|_{W^{s,p}(\Om,\CC)}\}$ is also a bounded sequence.  	
	\end{lemma}
	\begin{proof}
		\Cref{th-diamagnetic} implies that $\{[|f_n|]_{W^{s,p}(\Om,\CC)}\}_n$ is bounded. \Cref{std-poincare}, with the choice $g\equiv\frac{1}{\mathcal{L}^N(\Om)}$, and H\"older's inequality gives
		\begin{align*}
		\|f_n\|_{L^{p}(\Om,\CC)}
		\leq \||f_n|-\I{\Om}|f_n|g\|_{L^{p}(\Om,\CC)}+\|\I{\Om}|f_n|g\|_{L^{p}(\Om,\CC)}
		\leq C [|f_n|]_{W^{s,p}(\Om,\CC)}+C\|f_n\|_{L^q(\Om,\CC)}.
		\end{align*}
		Therefore, the sequence $\{ \|f_n\|_{L^p(\Om,\CC)} \}_n$ is bounded. \Cref{equality-magnetic-usual}, then, implies that $\{[f_n]_{W^{s,p}(\Om,\CC)}\}$ is bounded and consequently $\{ \|f_n\|_{W^{s,p}(\Om,\CC)} \}$ is bounded.	
	\end{proof}

	\begin{proof}[\textbf{Proof of \Cref{th-poincare-for-vf}}]
		We prove this result by method of contradiction. Suppose the statement of the theorem is false. Then there exist choices of $p,q,s,\delta$ and $\Om$, satisfying the hypotheses of the theorem, such that there is a sequence of functions $f_n\in W^{s,p}(\Om,\CC)$ such that  $d^q_A(f_n)\geq\delta\|f_n\|_{L^q(\Om,\CC)}$ but $\lim\limits_{n\to\infty}[f_n]_{W^{s,p}_A(\Om,\CC)}=E^{p,q}_{s,A}\|f_n\|_{L^q(\Om,\CC)}$. We normalize the sequence by considering $\{\frac{f_n}{\|f_n\|_{L^q(\Om,\CC)}}\}_n$ and still call it $f_n$, so that now it satisfies $\|f_n\|_{L^q(\Om,\CC)}=1$, $d^q_A(f_n)\geq\delta$ and $\lim\limits_{n\to\infty}[f_n]_{W^{s,p}_A(\Om,\CC)}=E^{p,q}_{s,A}$. Now is $\{[f_n]_{W^{s,p}_A(\Om,\CC)}\}$ is bounded. \Cref{seminorm-bdd-norm} implies that $\{ f_n \}$ is actually bounded in $W^{s,p}(\Om,\CC)$ with respect to its full norm. Hence, by \Cref{cpt-embedding}, there is some $f\in L^q(\Om,\CC)$ and a subsequence of $\{ f_n \}$ (without loss of generality, we assume it to be $\{ f_n \}$ itself) such that $f_n\to f$ strongly in $L^q(\Om,\CC)$, giving $\|f\|_{L^q(\Om,\CC)}=1$ and $d^q_{s,A}(f)\geq\delta$. Passing to a subsequence, we can assume $f_n\to f$ pointwise a.e. Finally using the fact $\lim\limits_{n\to\infty}[f_n]_{W^{s,p}_A(\Om,\CC)}=E^{p,q}_{s,A}$, Fatou's lemma and the definition of $E^{p,q}_{s,A}$ we get
		\begin{equation*}
		E^{p,q}_{s,A}=\lim\limits_{n\to \infty}[f_n]_{W^{s,p}_A(\Om,\CC)}\geq [f]_{W^{s,p}_A(\Om,\CC)} \geq E^{p,q}_{s,A}\|f\|_{L^q(\Om,\CC)}=E^{p,q}_{s,A}.
		\end{equation*}
		This shows that $f\in M^{p,q}_{s,A}$-which is a contradiction to the fact $d^q_{s,A}(f)\geq\delta$ for all $n$. Hence the first part of the theorem is proved.\smallskip
		
		Now we show that the best constant $S$ in \eqref{th-poincare-for-vf-eq1} is achieved. First note that \eqref{th-poincare-for-vf-eq1} can be rewritten as		
		\begin{equation}\label{th-poincare-for-vf-eq2}
		\left( \frac{1}{S}+E^{p,q}_{s,A} \right)\leq  \frac{[f]_{W^{s,p}_A(\Om,\CC)}}{\| f \|_{L^q(\Om,\CC)}}.
		\end{equation}
	If $S$ is the best possible constant, by definition, we have
		\begin{equation*}
			\frac{1}{S}=\inf_{\stackrel{f\in W^{s,p}_A(\Om,\CC)}{d^q_{s,A}(f)\geq \delta\| f \|_{L^q(\Om,\CC)}}}  \frac{[f]_{W^{s,p}_A(\Om,\CC)}}{\| f \|_{L^q(\Om,\CC)}}-E^{p,q}_{s,A}>0.
		\end{equation*} 
		So achievement of best constant $S$ in \eqref{th-poincare-for-vf-eq1} is equivalent to showing that the infimum (say $\xi$), taken over all $f\in W^{s,p}_A(\Om,\CC)$ with $d^q_{s,A}(f)\geq \delta\| f \|_{L^q(\Om,\CC)}$,  of RHS of \eqref{th-poincare-for-vf-eq2} is minimum. If this is not the case, then there exist a sequence $\{g_n\}$ in $W^{s,p}_A(\Om,\CC)$, with $\| g_n \|_{L^q(\Om,\CC)}=1$ and $d^q_{s,A}(g_n)\geq 1$, such that $[f]_{W^{s,p}_A(\Om,\CC)}\to\xi$. As in the previous part we can use the same argument to show that there exists $g\in W^{s,p}_A(\Om,\CC)$ such that $g_n\to g$ in $L^p(\Om,\CC)$, strongly together with $\| g \|_{L^q(\Om,\CC)}=1$. Again, as above, use of weak lower semicontinuity gives
		\begin{equation*}
		\xi=\lim\limits_{n\to \infty}[g_n]_{W^{s,p}_A(\Om,\CC)}\geq [g]_{W^{s,p}_A(\Om,\CC)} .
		\end{equation*}
		By definition of $\xi$, the above inequality is actually an equality and the proof concludes.
		
		To see the last part of the theorem, let, if possible, $M^{p,q}_{s,A}$ be empty. Then $d^q_{s,A}(f)=+\infty$, as it is the infimum of an empty set. Hence the condition $d^q_{s,A}(f)\geq \delta\| f \|_{L^q(\Om,\CC)}$ holds for any $\delta>0$ and for any $f\in L^p(\Om,\CC)$. We can then apply \eqref{th-poincare-for-vf-eq1} to any positive constant function (note that they all have the same $W^{s,p}_A$-seminorm). This contradicts the strict positivity of $S$.
	\end{proof}
	
	\begin{lemma}\label{poin-small-supp}
		Let $\Om,p,q$ be as in \Cref{std-poincare} and let $0<\delta_0\leq\delta<1$. Then there is a constant $S_1=S_1(\Om,\delta_0,p,q)>0$ such that $\forall\ f \in W^{s,p}(\Om,\CC)$ with $\mathcal{L}^N(\{x\ \Big| \ f(x)\neq 0\})\leq \mathcal{L}^N(\Om)(1-\delta)$ we have
		\begin{equation}
		\| f \|_{L^q(\Om,\CC)}\leq S_1 [f]_{W^{s,p}(\Om,\CC)}.
		\end{equation}
	\end{lemma}
	\begin{proof}
		\Cref{std-poincare} and H\"older's inequality imply
		\begin{align*}
		\| f \|_{L^q(\Om,\CC)}
		\leq& \left\| f-\frac{1}{\mathcal{L}^N(\Om)}\I{\Om}f \right\|_{L^q(\Om,\CC)}+\left\| \frac{1}{\mathcal{L}^N(\Om)}\I{\Om}f \right\|_{L^q(\Om,\CC)}\\
		=& \left\| f-\frac{1}{\mathcal{L}^N(\Om)}\I{\Om}f \right\|_{L^q(\Om,\CC)}+\left( \mathcal{L}^N(\Om) \right)^{\frac{1}{q}-1}\left|\I{\Om}f \right|\\
		\leq& C[f]_{W^{s,p}(\Om,\CC)}
		 + \left( \mathcal{L}^N(\Om) \right)^{\frac{1}{q}-1} \| f \|_{L^q(\Om,\CC)}\left( \mathcal{L}^N(\{x\ \Big| \ f(x)\neq 0\}) \right)^{1-\frac{1}{q}}\\
		\leq& C[f]_{W^{s,p}(\Om,\CC)} + (1-\delta)^{1-\frac{1}{q}}\| f \|_{L^q(\Om,\CC)}.
		\end{align*}
		Hence
		\begin{equation*}
		\| f \|_{L^q(\Om,\CC)}\leq \frac{C}{1-(1-\delta)^{1-\frac{1}{q}}} [f]_{W^{s,p}(\Om,\CC)}\leq \frac{C}{1-(1-\delta_0)^{1-\frac{1}{q}}} [f]_{W^{s,p}(\Om,\CC)}. 
		\end{equation*}
		This proves the lemma.
	\end{proof}

	\begin{proof}[\textbf{Proof of \Cref{th-poincare-punctured}}]
		Let, if possible, the statement of the theorem not be true. Then $\exists\ \varepsilon>0$ such that there are three sequences: $C_n>0$ with $C_n\to\infty$, measurable subsets $\Lambda_n\subseteq\Om$ and $f_n\in W^{s
			,p}(\Om,\CC)$ with $\|f_n\|_{L^q(\Om,\CC)} =1$ and $d^q_{s,A}(f_n)\geq\delta$ but 
		\begin{align}\label{th-poincare-punctured-eq1}
		\nonumber
		\left(\ \II{\Lambda_n\times\Lambda_n}\frac{|f_n(x)-e^{i(x-y)\cdot A\left(\frac{x+y}{2}\right)}f_n(y)|^p}{|x-y|^{N+sp}}dxdy\right)^\frac{1}{p}
		&+C_n\left(\ \II{(\Om\times\Om)\setminus(\Lambda_n\times\Lambda_n)}\frac{|f_n(x)-e^{i(x-y)\cdot A\left(\frac{x+y}{2}\right)}f_n(y)|^r}{|x-y|^{N+sr}}dxdy\right)^\frac{1}{r}\\
		<&\left(\frac{1}{S+\varepsilon}+E^{p,q}_{s,A}\right).
		\end{align}
		From this equation, it immediately follows that 
		\begin{equation}\label{th-poincare-punctured-eq9}
		\II{(\Om\times\Om)\setminus(\Lambda_n\times\Lambda_n)}\frac{|f_n(x)-e^{i(x-y)\cdot A\left(\frac{x+y}{2}\right)}f_n(y)|^r}{|x-y|^{N+sr}}dxdy\to 0 \mbox{ as } n\to\infty,
		\end{equation} 
		and 
		\begin{equation}\label{th-poincare-punctured-eq10}
		\II{\Lambda_n\times\Lambda_n}\frac{|f_n(x)-e^{i(x-y)\cdot A\left(\frac{x+y}{2}\right)}f_n(y)|^p}{|x-y|^{N+sp}}dxdy \mbox{ is bounded as } n\to\infty.
		\end{equation}
		Since $\Om$ is bounded, we assume $\mbox{diam}(\Om)<R$. The proof of the theorem is divided into the following three steps.\smallskip

		\noindent\textbf{Step-1:} $\Lambda_n\times\Lambda_n$ and $(\Om\times\Om)\setminus(\Lambda_n\times\Lambda_n)$ can be replaced by some $\lambda_n\ (\subseteq\Om\times\Om)$ and $\gamma_n=(\Om\times\Om)\setminus\lambda_n$ respectively, such that
		\begin{equation}\label{th-poincare-punctured-eq2}
		\sum_{n=1}^{\infty}\mathcal{L}^{2N}(\gamma_n)<\infty
		\end{equation}
		but still analogues for \eqref{th-poincare-punctured-eq9}, \eqref{th-poincare-punctured-eq10} hold.\smallskip

		Let us consider the set
		\begin{multline*}
		\gamma_n:=\left\{ (x,y)\in(\Om\times\Om)\setminus(\Lambda_n\times\Lambda_n)\ \Big| \ \frac{|f_n(x)-e^{i(x-y)\cdot A\left(\frac{x+y}{2} \right)}f_n(y)|}{|x-y|^{\frac{N}{r}+s}}  \right. \\
		 \left. \geq\left(\ \II{(\Om\times\Om)\setminus(\Lambda_n\times\Lambda_n)}\frac{|f_n(x)-e^{i(x-y)\cdot A\left(\frac{x+y}{2}\right)}f_n(y)|^r}{|x-y|^{N+sr}}dxdy  \right)^\frac{1}{2r} \right\}
		\end{multline*}
		and $\lambda_n :=(\Om\times\Om) \setminus \gamma_n$. 
		Then,
		\begin{align*}
		\II{(\Om\times\Om)\setminus(\Lambda_n\times\Lambda_n)}\frac{|f_n(x)-e^{i(x-y)\cdot A\left(\frac{x+y}{2}\right)}f_n(y)|^r}{|x-y|^{N+sr}}&dxdy
		\geq
		\II{\gamma_n}\frac{|f_n(x)-e^{i(x-y)\cdot A\left(\frac{x+y}{2}\right) }f_n(y)|^r}{|x-y|^{N+sr}}dxdy\\
		\geq&  \mathcal{L}^{2N}(\gamma_n)\Big( \II{(\Om\times\Om)\setminus(\Lambda_n\times\Lambda_n)}\frac{|f_n(x)-e^{i(x-y)\cdot A\left(\frac{x+y}{2}\right)}f_n(y)|^r}{|x-y|^{N+sr}}dxdy  \Big)^\frac{1}{2}.
		\end{align*}
		This shows that $\mathcal{L}^{2N}(\gamma_n)\to0$ as $n\to\infty$. Passing to a subsequence we can, without loss of generality, assume that 
		\eqref{th-poincare-punctured-eq2} holds. Note that in the second term in the LHS of \eqref{th-poincare-punctured-eq1}, the domain is going to be replaced by a smaller one, hence we need not worry about it. Therefore, we turn our attention to the first term of the LHS and show that the extra integral which is needed to be added there does not disrupt the inequality. Since $\mbox{diam}(\Om)<R$,
		\begin{align*}
		\II{((\Om\times\Om)\setminus(\Lambda_n\times\Lambda_n))\setminus\gamma_n}&\frac{|f_n(x)-e^{i(x-y)\cdot A\left(\frac{x+y}{2}\right)}f_n(y)|^p}{|x-y|^{N+sp}}dxdy\\
		=&\II{((\Om\times\Om)\setminus(\Lambda_n\times\Lambda_n))\setminus\gamma_n}\frac{|f_n(x)-e^{i(x-y)\cdot A\left(\frac{x+y}{2}\right)}f_n(y)|^p}{R^{N+sp}\left(\frac{|x-y|}{R}\right)^{N+sp}}dxdy\\
		\leq&
		C\II{((\Om\times\Om)\setminus(\Lambda_n\times\Lambda_n))\setminus\gamma_n}\frac{|f_n(x)-e^{i(x-y)\cdot A\left(\frac{x+y}{2}\right)}f_n(y)|^p}{R^{N+sp}\left(\frac{|x-y|}{R}\right)^{\frac{Np}{r}+sp}}dxdy\\
		\leq& C(R,p,r,s)
		C\mathcal{L}^{2N}(\Om\times\Om)
		\Big( \II{(\Om\times\Om)\setminus(\Lambda_n\times\Lambda_n)}\frac{|f_n(x)-e^{i(x-y)\cdot A\left(\frac{x+y}{2}\right)}f_n(y)|^r}{|x-y|^{N+sr}}dxdy  \Big)^\frac{p}{2r}
		\to 0.
		\end{align*}
		Therefore, \eqref{th-poincare-punctured-eq1} can be rewritten as
		\begin{multline}\label{th-poincare-punctured-eq3}
		\left(\II{\lambda_n}\frac{|f_n(x)-e^{i(x-y)\cdot A\left(\frac{x+y}{2}\right)}f_n(y)|^p}{|x-y|^{N+sp}}dxdy\right)^\frac{1}{p}
		+C_n\left(\II{\gamma_n}\frac{|f_n(x)-e^{i(x-y)\cdot A\left(\frac{x+y}{2}\right)}f_n(y)|^r}{|x-y|^{N+sr}}dxdy\right)^\frac{1}{r}\\
		<\left(\frac{1}{S+\varepsilon}+E^{p,q}_{s,A}\right).
		\end{multline}
		This implies (similarly as \eqref{th-poincare-punctured-eq9} and\eqref{th-poincare-punctured-eq10})
		\begin{equation}\label{th-poincare-punctured-eq11}
		\II{\gamma_n}\frac{|f_n(x)-e^{i(x-y)\cdot A\left(\frac{x+y}{2}\right)}f_n(y)|^r}{|x-y|^{N+sr}}dxdy\to 0 \mbox{ as } n\to\infty
		\end{equation} 
		and 
		\begin{equation}\label{th-poincare-punctured-eq12}
		\II{\lambda_n}\frac{|f_n(x)-e^{i(x-y)\cdot A\left(\frac{x+y}{2}\right)}f_n(y)|^p}{|x-y|^{N+sp}}dxdy \mbox{ is bounded as } n\to\infty
		\end{equation} 
		as before.\smallskip

		\noindent\textbf{Step-2:} $f_n\to f$ strongly in $L^q(\Om,\CC)$.\smallskip
		
		
		\noindent\textbf{Case-I:} Assume that condition (a) or (b) in the hypotheses hold. That is $q<\frac{N}{N-s}$ (then we have $r=1$) or $q=\infty$ (then $sr\geq N$).
		\smallskip
		
		Since $r<p$, \eqref{th-poincare-punctured-eq11} and \eqref{th-poincare-punctured-eq12} imply that $[f_n]_{W^{s,r}_A(\Om,\CC)}$ is uniformly bounded. \Cref{seminorm-bdd-norm} implies that $\|f_n\|_{W^{s,r}(\Om,\CC)}$ is bounded. Thus we can apply \Cref{cpt-embedding} to complete the step.\smallskip

		\noindent\textbf{Case-II:} 
		Assume that condition (c) in the hypotheses hold. Then $q\in\left[\frac{N}{N-s},\infty\right)$,$sr<N$ and $q<\frac{Nr}{N-sr}$.\smallskip
		
		\noindent\textbf{Step-2/a:} $f_n\rightharpoonup f$ in $L^q(\Om,\CC)$ for some $f$.\smallskip
		
		Since the sequence $\{f_n\}$ is bounded in $L^q(\Om,\CC)$ and hence there is a subsequence, still denoted by
		$f_n$, and an $f \in L^q(\Om,\CC)$, such that $f_n\rightharpoonup f$ in $L^q(\Om,\CC)$.

		\label{originality}\noindent\textbf{Step-2/b:} There exists $\eta_0>0$ such that for any $\eta_0>\eta>0$, there exists suitable $\alpha>0$ such that $[f_n]_{W^{\alpha s,r}}(\Om,\CC)$ is bounded, hence, up to a subsequence, it converges to some $f\in L^{q-\eta}(\Om,\CC)$.\smallskip

		Note that $1\leq r<q<\frac{Nr}{N-sr}<\infty$. Choose $\eta>0$, small enough, such that $q-r-\eta>0$ which implies $q-\eta>1$. Fix $\alpha\in\left(\frac{N(q-r-\eta)}{sr(q-\eta)},\mbox{min}\left\{1,\frac{N}{sr}\right\}\right)$. This gives $\alpha s r<N$ and $1<q-\eta<\frac{Nr}{N-\alpha s r}$. Hence \Cref{cpt-embedding} gives
		\begin{equation}\label{th-poincare-punctured-eq13}
		W^{\alpha s,r}(\Om,\CC) \mbox{ is compactly embedded in } L^{q-\eta}(\Om,\CC).
		\end{equation} 
		Note that both the terms in the LHS of \eqref{th-poincare-punctured-eq3} are bounded. Since $\alpha\in(0,1)$ and $\alpha r< r\leq p$, applying \Cref{p-r-embedding}, we get that the terms
		$\left(\II{\lambda_n}\frac{|f_n(x)-e^{i(x-y)\cdot A\left(\frac{x+y}{2}\right)}f_n(y)|^r}{|x-y|^{N+\alpha sr}}dxdy\right)^\frac{1}{r}$ and 
		$\left(\II{\gamma_n}\frac{|f_n(x)-e^{i(x-y)\cdot A\left(\frac{x+y}{2}\right)}f_n(y)|^r}{|x-y|^{N+\alpha sr}}dxdy\right)^\frac{1}{r}$ are bounded. This implies that $\{f_n\}$ is bounded in $W^{\alpha s,r}(\Om,\CC)$. By \eqref{th-poincare-punctured-eq13}, $\{f_n\}$ has a strongly convergent subsequence in $L^{q-\eta}(\Om,\CC)$, which we assume to be itself, converging to $f$. (Note that this $f$ is actually what we talked about in step-2/a.)

		\noindent\textbf{Step-2/c:} Estimates on $q$ norm of $f_n$ to show that $f_n\to f$ in $L^q(\Om,\CC)$.\smallskip
		
		Take $M>1$ such that $M^{-q}<\mathcal{L}^N(\Om)(1-\frac{1}{2})$ and define, for $x\in\Om$,
		\begin{equation*}
		f^M_n(x)=\min\{M,|f_n(x)|\} \quad \mbox{and} \quad  h^M_n(x)=|f_n(x)|-f^M_n(x)\geq 0.
		\end{equation*}
		Now observe that 
		\begin{equation*}
		|f_n(x)|>M \quad \mbox{if and only if}\quad h^M_n(x)\neq 0.
		\end{equation*} 
		We use this fact to derive,
		\begin{align}\label{th-poincare-punctured-eq6}
		\mathcal{L}^N\left(|f_n|^{-1}((M,\infty))\right)
		=&\mathcal{L}^N(\{ x\in\Om\ \Big|\ h^M_n(x)\neq 0 \})\\
		\nonumber
		=& \I{h_n^M(x)\neq 0}dx
		\leq \frac{\|f_n\|_{L^q(\Om,\CC)}^q}{M^q}	
		\leq M^{-q}<\mathcal{L}^N(\Om)(1-\frac{1}{2}).
		\end{align}
		Choose suitable $s_1\in(0,s),\ r_1\in[1,r)$ such that $s_1r_1<N$ and $q=	\frac{Nr_1}{N-s_1r_1}$. We apply \Cref{poin-small-supp} on $W^{s_1,r_1}(\Om,\CC)$ and $L^{q}(\Om,\CC)$. So,
		\begin{equation}\label{th-poincare-punctured-eq7}
		\|h^M_n\|_{L^{q}(\Om,\CC)}\leq S_1(\Om,N,s,q,r)[h^M_n]_{W^{s_1,r_1}(\Om,\CC)}.
		\end{equation}
		Define $H^M_n:=(\Om\times\Om)\setminus \{ (x,y)\in\Om\times\Om\ \Big|\ |f_n(x)|,|f_n(y)|\leq M \}$. 
		We use symmetry of the integrand and  \Cref{th-diamagnetic} in the following computation 
		\begin{align*}
		[h^M_n]_{W^{s_1,r_1}(\Om,\CC)}^{r_1}
		=&\II{\Om\times\Om}\frac{|h^M_n(x)-h^M_n(y)|^{r_1}}{|x-y|^{N+s_1r_1}}dxdy\\
		=&\II{\substack{\Om\times\Om\\|f_n(x)|,\ |f_n(y)|\geq M}}\frac{|h^M_n(x)-h^M_n(y)|^{r_1}}{|x-y|^{N+s_1r_1}}dxdy
		+2 \II{\substack{\Om\times\Om\\|f_n(x)|\geq M>|f_n(y)|}}\frac{|h^M_n(x)-h^M_n(y)|^{r_1}}{|x-y|^{N+s_1r_1}}dxdy\\
		& +\II{\Om\times\Om \setminus H^M_n}\frac{|h^M_n(x)-h^M_n(y)|^{r_1}}{|x-y|^{N+s_1r_1}}dxdy\\
		=&\II{\substack{\Om\times\Om\\|f_n(x)|,\ |f_n(y)|\geq M}}\frac{||f_n(x)|-|f_n(y)||^{r_1}}{|x-y|^{N+s_1r_1}}dxdy
		+2\II{\substack{\Om\times\Om\\|f_n(x)|\geq M>|f_n(y)|}}\frac{||f_n(x)|-M|^{r_1}}{|x-y|^{N+s_1r_1}}dxdy\\
		\leq& \II{\substack{\Om\times\Om\\|f_n(x)|,\ |f_n(y)|\geq M}}\frac{||f_n(x)|-|f_n(y)||^{r_1}}{|x-y|^{N+s_1r_1}}dxdy
		+2\II{\substack{\Om\times\Om\\|f_n(x)|\geq M>|f_n(y)|}}\frac{||f_n(x)|-|f_n(y)||^{r_1}}{|x-y|^{N+s_1r_1}}dxdy\\
		=& \II{H^M_n}\frac{||f_n(x)|-|f_n(y)||^{r_1}}{|x-y|^{N+s_1r_1}}dxdy\\
		\leq& \II{H^M_n}\frac{|f_n(x)-e^{i(x-y)\cdot A\left(\frac{x+y}{2}\right)}f_n(y)|^{r_1}}{|x-y|^{N+s_1r_1}}dxdy\\
		\leq& \II{H^M_n\cap\lambda_n}\frac{|f_n(x)-e^{i(x-y)\cdot A\left(\frac{x+y}{2}\right)}f_n(y)|^{r_1}}{|x-y|^{N+s_1r_1}}dxdy
		+ \II{H^M_n\cap\gamma_n}\frac{|f_n(x)-e^{i(x-y)\cdot A\left(\frac{x+y}{2}\right)}f_n(y)|^{r_1}}{|x-y|^{N+s_1r_1}}dxdy.
		\end{align*}
		Since $r_1<r$, \Cref{p-r-embedding} and \eqref{th-poincare-punctured-eq11} imply that the last term converges to $0$. Thus the above calculation and \eqref{th-poincare-punctured-eq7} then imply
		\begin{align*}
		&\limsup\limits_{n\to \infty}\|h^M_n\|_{L^{q}(\Om,\CC)}^{r_1}\\
		\leq& S_1^{r_1}\limsup\limits_{n\to \infty}[h^M_n]_{W^{s_1,r_1}(\Om,\CC)}^{r_1}\\
		\leq& S_1^{r_1}\limsup\limits_{n\to \infty} \II{H^M_n\cap\lambda_n}\frac{|f_n(x)-e^{i(x-y)\cdot A\left(\frac{x+y}{2}\right)}f_n(y)|^{r_1}}{|x-y|^{N+s_1r_1}}dxdy\\
		\leq& S_1^{r_1}\limsup\limits_{n\to \infty} \left(\II{\lambda_n}\frac{|f_n(x)-e^{i(x-y)\cdot A\left(\frac{x+y}{2}\right)}f_n(y)|^p}{|x-y|^{N+sp}}dxdy\right)^\frac{r_1}{p} \left(\II{H^M_n}|x-y|^{-N+\frac{sr_1p-s_1r_1p}{p-r_1}}dxdy \right)^\frac{p-r_1}{p}\\
		\leq& S_1^{r_1}\left(\frac{1}{S+\varepsilon}+E^{p,q}_{s,A}\right)^{r_1}
		\limsup\limits_{n\to \infty} \left(\II{H^M_n}|x-y|^{-N+\frac{sr_1p-s_1r_1p}{p-r_1}}dxdy \right)^\frac{p-r_1}{p}.
		\end{align*}
		It is clear that $H^M_n$ is the disjoint union of following three sets:\\ $\left(|f_n|^{-1}([M,\infty))\times |f_n|^{-1}([M,\infty))\right)$, $\left(|f_n|^{-1}([M,\infty))\times |f_n|^{-1}([0,M))\right)$ and \\ $\left(|f_n|^{-1}([0,M))\times |f_n|^{-1}([M,\infty))\right)$. Along with this decomposition of $H^M_n$, we use \eqref{th-poincare-punctured-eq6} to get
		\begin{align*}
		\II{H^M_n}|x-y|^{-N+\frac{sr_1p-s_1r_1p}{p-r_1}}dxdy
		&\leq \II{|f_n|^{-1}([M,\infty))\times B(x,R)}|x-y|^{-N+\frac{sr_1p-s_1r_1p}{p-r_1}}dxdy\\
		&\ +2\II{|f_n|^{-1}([M,\infty))\times B(x,R)}|x-y|^{-N+\frac{sr_1p-s_1r_1p}{p-r_1}}dxdy\\
		&\leq C(N,p,r,s,R)M^{-q}.
		\end{align*}
	In the above calculation, we have used the facts that $\mbox{diam}(\Om)<R$ and $s_1<s$. Now we have
		\begin{equation}\label{th-poincare-punctured-eq8}
		\limsup\limits_{n\to \infty}\|h^M_n\|_{L^{q}(\Om,\CC)}\leq S_1\left(\frac{1}{S+\varepsilon}+E^{p,q}_{s,A}\right)M^\frac{q(r_1-p)}{pr_1}.
		\end{equation}
		Therefore \eqref{th-poincare-punctured-eq1} and the triangle inequality of $L^q$-norm gives
		\begin{align*}
		\I{\Om}|f|^{q-\eta}
		=&\lim\limits_{n\to\infty}\I{\Om}|f_n|^{q-\eta}
		\geq \lim\limits_{n\to\infty}\I{\Om}|f^M_n|^{q-\eta}
		\geq M^{-\eta}\lim\limits_{n\to\infty}\I{\Om}|f^M_n|^q\\
		\geq& M^{-\eta}\lim\limits_{n\to\infty}\left[ \|f_n\|_{L^q(\Om,\CC)}-\|h^M_n\|_{L^q(\Om,\CC)} \right]^q\\
		\geq& M^{-\eta}\left( 1- S_1\left(\frac{1}{S+\varepsilon}+E^{p,q}_{s,A}\right)M^\frac{q(r_1-p)}{pr_1}
		\right)^q,
		\end{align*}
		which after taking the limit as $\eta\to 0$ and then $M\to\infty$ gives (recall that $r_1<p$)
		$$
		\I{\Om}|f|^q\geq 1 \mbox{ and hence } \I{\Om}|f|^q = 1.
		$$
		So, $f_n\rightharpoonup f$ and $\|f_n\|_{L^q(\Om,\CC)}\to \|f\|_{L^q(\Om,\CC)}$, implying the strong convergence $f_n\to f$ in $L^q(\Om,\CC)$.
		\smallskip

		\noindent\textbf{Step-3:} The strong convergence $f_n\to f$ in $L^q(\Om,\CC)$ gives a contradiction.\smallskip

		Now the fact that $f_n\to f$ in $L^q(\Om,\CC)$ implies that  $f_n(x)\to f(x)$ for almost all $x\in\Om$. For fixed $k\in\NN$, set $\Sigma_k:=(\Om\times\Om)\setminus (\gamma_k\cup\gamma_{k+1}\cup\cdots)$. Then for $n\geq k $, we have, using Fatou's lemma in the last inequality,
		\begin{align*}
		\liminf\limits_{n\to \infty}\II{\lambda_n}\frac{|f_n(x)-e^{i(x-y)\cdot A\left(\frac{x+y}{2}\right)}f_n(y)|^p}{|x-y|^{N+sp}}dxdy
		=&\liminf\limits_{n\to \infty}\II{(\Om\times\Om)\setminus\gamma_n}\frac{|f_n(x)-e^{i(x-y)\cdot A\left(\frac{x+y}{2}\right)}f_n(y)|^p}{|x-y|^{N+sp}}dxdy\\
		\geq& 
		\liminf\limits_{n\to\infty} \II{\Sigma_k}\frac{|f_n(x)-e^{i(x-y)\cdot A\left(\frac{x+y}{2}\right)}f_n(y)|^p}{|x-y|^{N+sp}}dxdy\\
		\geq& \II{\Sigma_k}\frac{|f(x)-e^{i(x-y)\cdot A\left(\frac{x+y}{2}\right)}f(y)|^p}{|x-y|^{N+sp}}dxdy.
		\end{align*}
		This holds for any $k\in\NN$. Also $\Sigma_k\subseteq\Sigma_{k+1}\subseteq\Om\times\Om$ and $\mathcal{L}^{2N}(\bigcup_{k=0}^{\infty}\Sigma_k)=\mathcal{L}^{2N}(\Om\times\Om)$ by \eqref{th-poincare-punctured-eq2}. Therefore
		\begin{equation}\label{th-poincare-punctured-eq5}
		\liminf\limits_{n\to \infty}\II{\lambda_n}\frac{|f_n(x)-e^{i(x-y)\cdot A\left(\frac{x+y}{2}\right)}f_n(y)|^p}{|x-y|^{N+sp}}dxdy
		\geq \II{\Om\times\Om}\frac{|f(x)-e^{i(x-y)\cdot A\left(\frac{x+y}{2}\right)}f(y)|^p}{|x-y|^{N+sp}}dxdy.
		\end{equation}
		Since $f_n\to f$ in $L^q(\Om,\CC)$, we have $d^q_{s,A}(f)\geq\delta\|f\|_{L^q(\Om,\CC)}=\delta$. Now \Cref{th-poincare-for-vf} and \eqref{th-poincare-punctured-eq5} imply
		\begin{equation*}
		\liminf\limits_{n\to \infty}\left(\ \II{\lambda_n}\frac{|f_n(x)-e^{i(x-y)\cdot A\left(\frac{x+y}{2}\right)}f_n(y)|^p}{|x-y|^{N+sp}}dxdy\right)^\frac{1}{p}\geq  \left(\frac{1}{S}+E^{p,q}_{s,A}\right)
		\end{equation*}
		which contradicts \eqref{th-poincare-punctured-eq3} and  hence the proof follows.
	\end{proof}
	\bigskip

	\section{The magnetic fractional Laplacian operator}\label{laplacian}
	In this section, we consider the case $p=2$. $A$ shall, denote a bounded, vector field on the convex hull of the domain $\Om$. We first observe the following:
	
	\begin{lemma}\label{self adjoint}
		$(-\Delta_{2,A})^s$ is a self adjoint operator on $L^2(\Om,\CC)$.
	\end{lemma}
	\begin{proof} The result follows from the following calculation:
		\begin{align*}
		<(-\Delta_{p,A})^sf,g>
		=&\II{\Om\times\Om} \frac{(f(x)-e^{i(x-y)\cdot A\left(\frac{x+y}{2}\right)}f(y))}{|x-y|^{N+2s}}\overline{g(x)}dydx\\
		=&\II{\Om\times\Om} \frac{(f(x)-e^{i(x-y)\cdot A\left(\frac{x+y}{2}\right)}f(y))}{2|x-y|^{N+2s}}\overline{g(x)}dydx
		+\II{\Om\times\Om} \frac{(f(y)-e^{i(y-x)\cdot A\left(\frac{x+y}{2}\right)}f(x))}{2|x-y|^{N+2s}}\overline{g(y)}dydx\\
		=&\II{\Om\times\Om} \frac{(f(x)-e^{i(x-y)\cdot A\left(\frac{x+y}{2}\right)}f(y))(\overline{g(x)}-e^{i(y-x)\cdot A\left(\frac{x+y}{2}\right)}\overline{g(y)})}{2|x-y|^{N+2s}}dydx\\
		=&\II{\Om\times\Om} \frac{(\overline{g(x)}-e^{i(y-x)\cdot A\left(\frac{x+y}{2}\right)}\overline{g(y)})}{2|x-y|^{N+2s}}f(x)dydx
		-\II{\Om\times\Om} \frac{(\overline{g(x)}e^{i(x-y)\cdot A\left(\frac{x+y}{2}\right)}-\overline{g(y)})}{2|x-y|^{N+2s}}f(y)dydx\\
		=&\II{\Om\times\Om} \frac{(\overline{g(x)}-e^{i(y-x)\cdot A\left(\frac{x+y}{2}\right)}\overline{g(y)})}{|x-y|^{N+2s}}f(x)dydx\\
		=&\I{\Om}f(x)\overline{(-\Delta_{p,A})^sg(x)}dx=<f,(-\Delta_{p,A})^sg>.
		\end{align*}
		
	\end{proof}
\smallskip

	We define $E_1:=\sqrt{E^{2,2}_{s,A}}$ to be the first eigenvalue of $(-\Delta_{2,A})^s$ and we denote the corresponding eigenfunction by $\phi_1$. Since the eigenvectors of a self-adjoint operator, on a Hilbert space, form an orthonormal basis, consecutive eigenpairs $(E_n,\phi_n)$ of $(-\Delta_{2,A})^s$ can be defined iteratively, as follows: if $(E_n,\phi_n)$ is the $n$'th eigenpair, we define
	$$
	E_{n+1}:=\inf \left\{ [f]_{W^{s,2}_A(\Om,\CC)}^2\ \Big| \ \|f\|_{L^2(\Om,\CC)}=1,\ \I{\Om}f\overline{\phi_j}=0 \mbox{ for all } 1\leq j\leq n \right\}.
	$$
	\Cref{cpt-embedding} ensures that a minimizer of $E_{n+1}$, which we call $\phi_{n+1}$, must exist.
	\smallskip
	
	\begin{proposition}\label{discrete}
		Let $\Om\subseteq \RR^N$ be a bounded Lipschitz domain. The set of eigenvalues of $(-\Delta_{2,A})^s$ in $W^{s,2}(\Om,\CC)$ is discrete.
	\end{proposition}
	
	\begin{proof}
		If possible, let the statement not be true. Then there exists some $E\in(0,\infty)$ such that there are infinitely many (all the) eigenvalues of $(-\Delta_{2,A})^s$ below $E$. We consider the sequence of eigenfunctions $\{\phi_n\}_n$ in $W^{s,2}(\Om,\CC)$. Now $\phi_n$ are all of unit $L^2$-norm and that they achieve the corresponding eigenvalues. We use this fact in the following calculation, where the first line is achieved by adding and subtracting the term $\phi_n(y)e^{iA(\frac{x+y}{2}\cdot (x-y))}$ in the numerator of the integrand.
		\begin{align*}
		[\phi_n]_{W^{s,2}(\Om,\CC)}^2
		=&\II{\Om\times\Om}\frac{|\phi_n(x)-\phi_n(y)|^2}{|x-y|^{N+2s}}dxdy\\
		\leq& C\II{\Om\times\Om}\frac{|\phi_n(x)-\phi_n(y)e^{iA(\frac{x+y}{2}\cdot (x-y))}|^2}{|x-y|^{N+2s}}dxdy
		+C\|\phi_n\|_{L^\infty(\Om)}\II{\Om\times\Om}\frac{|e^{iA(\frac{x+y}{2}\cdot (x-y))}-1|^2}{|x-y|^{N+2s}}dxdy\\
		=&C[\phi_n]_{W^{s,2}_A(\Om,\CC)}^2+C[1]_{W^{s,2}_A(\Om,\CC)}^2
		\leq CE+C[1]_{W^{s,2}_A(\Om,\CC)}^2.
		\end{align*}
		This implies that $\{\phi_n\}_n$ is a bounded sequence in $W^{s,2}(\Om,\CC)$. By \Cref{cpt-embedding}, we then, conclude that up to a subsequence, the sequence of orthonormal functions, $\{\phi_n\}_n$, converges in $L^2(\Om,\CC)$. This indeed is a contradiction and hence the result follows.
	\end{proof}

	\bigskip

	\section*{Funding:} \emph{
	Research work of the first author is funded by Matrics grant  (MTR/2020/000594).\\
	Research work of the second author is funded by Academy of Finland (Suomen Akatemia) grant: Geometrinen Analyysi(21000046081).\\
	Research work of the third author is funded by Matrics grant: (MTR/2019/000585) and by Core Research Grant (CRG/2022/007867) of SERB.}
	\smallskip

	\section*{Conflict of Interest}
	\emph{The authors have no competing interests to declare that are relevant to the content of this article.}
	\vfill\pagebreak

	\appendix
	\section{} \label{appendix1}
	
	\begin{lemma}
		For the sequence of functions, $f_\varepsilon$, defined in \eqref{example-function}, $$[f_\varepsilon]_{W^{s,r}(\Gamma,\CC)}\to 0 \quad \mbox{as } \varepsilon\to0.$$
	\end{lemma}
	\begin{proof}
		Denoting $\Gamma_1:=(0,\varepsilon]\times(0,1)$ and $\Gamma_2:=(\varepsilon,1)\times(0,1)$, we have 
		\begin{align*}
		[f_\varepsilon]_{W^{s,r}(\Gamma,\CC)}^r
		=&\II{\Gamma_1\times\Gamma_1}\frac{|f_\varepsilon(x)-f_\varepsilon(y)|^r}{|x-y|^{2+sr}}dxdy
		+ 2\II{\Gamma_1\times\Gamma_2}\frac{|f_\varepsilon(x)-f_\varepsilon(y)|^r}{|x-y|^{2+sr}}dxdy\\
		=&2\II{\Gamma_1\times\Gamma}\frac{|f_\varepsilon(x)-f_\varepsilon(y)|^r}{|x-y|^{2+sr}}dxdy\\
		=&2\I{x\in\Gamma_1}\I{\stackrel{y\in\Gamma}{|x-y|<\varepsilon}}\frac{|f_\varepsilon(x)-f_\varepsilon(y)|^r}{|x-y|^{2+sr}}dydx
		+2\I{x\in\Gamma_1}\I{\stackrel{y\in\Gamma}{|x-y|\geq\varepsilon}}\frac{|f_\varepsilon(x)-f_\varepsilon(y)|^r}{|x-y|^{2+sr}}dydx\\
		=&2\I{x\in\Gamma_1}\left(I_1(x)+I_2(x)\right)dx.
		\end{align*}
		Note that $f_\varepsilon$ is a Lipschitz function with Lipschitz constant $\frac{2(2-\varepsilon)}{\varepsilon(2+\varepsilon)}$. Also $|f_\varepsilon|\leq 1$ on $\Gamma$. We are going to use these facts in the following two estimates.
		
		\begin{equation*}
		I_1(x)
		=\I{\stackrel{y\in\Gamma}{|x-y|<\varepsilon}}\frac{|f_\varepsilon(x)-f_\varepsilon(y)|^r}{|x-y|^{2+sr}}dy
		\leq C\varepsilon^{-r} \I{\stackrel{y\in\Gamma}{|x-y|<\varepsilon}}\frac{dy}{|x-y|^{2+sr-r}}
		\leq \frac{C}{\varepsilon^{r}}\I{t=0}^\varepsilon t^{r-sr-1}dt=C\varepsilon^{-sr}.
		\end{equation*}
		
		\begin{equation*}
		I_2(x)
		=\I{\stackrel{y\in\Gamma}{|x-y|\geq \varepsilon}}\frac{|f_\varepsilon(x)-f_\varepsilon(y)|^r}{|x-y|^{2+sr}}dy
		\leq \I{\stackrel{y\in\Gamma}{|x-y|\geq\varepsilon}}\frac{dy}{|x-y|^{2+sr}}
		\leq C\I{t=\varepsilon}^2 t^{-sr-1}dt=C(\varepsilon^{-sr}-2^{-sr}).
		\end{equation*}
		Combining the above three calculations, we get $$[f_\varepsilon]_{W^{s,r}(\Gamma,\CC)}^r\leq 2C(2\varepsilon^{1-sr}-2^{-sr}\varepsilon)\to 0 \quad \mbox{as } \varepsilon\to 0.$$
	\end{proof}
\smallskip 

\begin{lemma}\label{ap-char}
Let $\Om=B(0,1)$, $\Gamma=B(0,1)\setminus B(0,\frac{1}{2})$, $\Lambda= B(0,\frac{1}{2})$ and	$f=\chi_{\Lambda}$. Then $f\in W^{s,p}(\Om)$, provided $sp<1$.
\end{lemma}
\begin{proof}
	The proof follows from the following calculation:
	\begin{align*}
		\frac{1}{2}[f]_{W^{s,p}(\Om)}^p
		=&\I{x\in\Lambda}\I{y\in\Gamma}\frac{dy}{|x-y|^{N+sp}}dx
		=\I{x\in\Lambda}\I{y\in-x+\Gamma}\frac{dy}{|y|^{N+sp}}dx
		\leq C(N)\I{x\in\Lambda}\I{r=\frac{1}{2}-|x|}^2\frac{r^{N-1}dr}{r^{N+sp}}dx\\
		=& C(N,s,p)\I{x\in B(0,\frac{1}{2})}\left(\left(\frac{1}{2}-|x|\right)^{-sp}-2^{-sp}\right)dx
		\leq C(N,s,p)\I{r=0}^\frac{1}{2}\left(\frac{1}{2}-r\right)^{-sp}r^{N-1}dr\\
		\leq&  C(N,s,p)\I{r=0}^\frac{1}{2}\left(\frac{1}{2}-r\right)^{-sp}dr
		\leq  C(N,s,p)\I{r=0}^\frac{1}{2}r^{-sp}dr
		<\infty,
	\end{align*}
where the last inequality follows from the hypothesis $sp<1$.
\end{proof}

	\bigskip

	\bibliography{bibliography}

\begin{thebibliography}{10}

\bibitem{AbFaTe}
Nicola Abatangelo, Mouhamed~Moustapha Fall, and Remi~Yvant Temgoua.
\newblock A {H}opf lemma for the regional fractional {L}aplacian.
\newblock {\em Ann. Mat. Pura Appl. (4)}, 202(1):95--113, 2023.

\bibitem{Ambrosio22}
Vincenzo Ambrosio.
\newblock Concentration phenomena for fractional magnetic {NLS} equations.
\newblock {\em Proc. Roy. Soc. Edinburgh Sect. A}, 152(2):479--517, 2022.

\bibitem{AmdAv}
Vincenzo Ambrosio and Pietro d'Avenia.
\newblock Nonlinear fractional magnetic {S}chr\"{o}dinger equation: existence
  and multiplicity.
\newblock {\em J. Differential Equations}, 264(5):3336--3368, 2018.

\bibitem{AnWa}
Harbir Antil and Mahamadi Warma.
\newblock Optimal control of the coefficient for the regional fractional
  {$p$}-{L}aplace equation: approximation and convergence.
\newblock {\em Math. Control Relat. Fields}, 9(1):1--38, 2019.

\bibitem{BMRS}
Kaushik Bal, Kaushik Mohanta, Prosenjit Roy, and Firoj Sk.
\newblock Hardy and {P}oincar\'e inequalities in fractional {O}rlicz-{S}obolev
  space.
\newblock {\em arXiv preprint arxiv:2009.07035}, 2020.

\bibitem{BBM}
Jean Bourgain, Haim Brezis, and Petru Mironescu.
\newblock Another look at {S}obolev spaces.
\newblock In {\em Optimal control and partial differential equations}, pages
  439--455. IOS, Amsterdam, 2001.

\bibitem{Bre11}
Haim Brezis.
\newblock {\em Functional analysis, {S}obolev spaces and partial differential
  equations}.
\newblock Universitext. Springer, New York, 2011.

\bibitem{Ind}
Indranil Chowdhury, Gyula Csat\'{o}, Prosenjit Roy, and Firoj Sk.
\newblock Study of fractional {P}oincar\'{e} inequalities on unbounded domains.
\newblock {\em Discrete Contin. Dyn. Syst.}, 41(6):2993--3020, 2021.

\bibitem{dASq}
Pietro d'Avenia and Marco Squassina.
\newblock Ground states for fractional magnetic operators.
\newblock {\em ESAIM Control Optim. Calc. Var.}, 24(1):1--24, 2018.

\bibitem{hhg}
Eleonora Di~Nezza, Giampiero Palatucci, and Enrico Valdinoci.
\newblock Hitchhiker's guide to the fractional {S}obolev spaces.
\newblock {\em Bull. Sci. Math.}, 136(5):521--573, 2012.

\bibitem{evans}
Lawrence~C. Evans.
\newblock {\em Partial differential equations}, volume~19 of {\em Graduate
  Studies in Mathematics}.
\newblock American Mathematical Society, Providence, RI, second edition, 2010.

\bibitem{Fall}
Mouhamed~Moustapha Fall.
\newblock Regional fractional {L}aplacians: boundary regularity.
\newblock {\em J. Differential Equations}, 320:598--658, 2022.

\bibitem{BoSa21}
Juli\'{a}n Fern\'{a}ndez~Bonder and Ariel~M. Salort.
\newblock Magnetic fractional order {O}rlicz-{S}obolev spaces.
\newblock {\em Studia Math.}, 259(1):1--24, 2021.

\bibitem{Guan}
Qing-Yang Guan.
\newblock Integration by parts formula for regional fractional {L}aplacian.
\newblock {\em Comm. Math. Phys.}, 266(2):289--329, 2006.

\bibitem{GuMe}
Z.~Guo and M.~Melgaard.
\newblock Fractional magnetic {S}obolev inequalities with two variables.
\newblock {\em Math. Inequal. Appl.}, 22(2):703--718, 2019.

\bibitem{HuVa}
Ritva Hurri-Syrj\"{a}nen and Antti~V. V\"{a}h\"{a}kangas.
\newblock On fractional {P}oincar\'{e} inequalities.
\newblock {\em J. Anal. Math.}, 120:85--104, 2013.

\bibitem{Kum}
S.~Kumaresan.
\newblock {\em Topology of metric spaces}.
\newblock Narosa Publishing House, New Delhi, 2005.

\bibitem{LiLo}
Elliott~H. Lieb and Michael Loss.
\newblock {\em Analysis}, volume~14 of {\em Graduate Studies in Mathematics}.
\newblock American Mathematical Society, Providence, RI, second edition, 2001.

\bibitem{LiSe}
Elliott~H Lieb and Robert Seiringer.
\newblock Proof of bose-einstein condensation for dilute trapped gases.
\newblock {\em Physical review letters}, 88(17):170409, 2002.

\bibitem{LiSeSoYn}
Elliott~H. Lieb, Robert Seiringer, Jan~Philip Solovej, and Jakob Yngvason.
\newblock {\em The mathematics of the {B}ose gas and its condensation},
  volume~34 of {\em Oberwolfach Seminars}.
\newblock Birkh\"{a}user Verlag, Basel, 2005.

\bibitem{LiSeYn02}
Elliott~H Lieb, Robert Seiringer, and Jakob Yngvason.
\newblock Superfluidity in dilute trapped bose gases.
\newblock In {\em The Stability of Matter: From Atoms to Stars}, pages
  903--908. Springer, 2002.

\bibitem{LiSeYn}
Elliott~H. Lieb, Robert Seiringer, and Jakob Yngvason.
\newblock Poincar\'{e} inequalities in punctured domains.
\newblock {\em Ann. of Math. (2)}, 158(3):1067--1080, 2003.

\bibitem{LiChGu}
Min Liu, Deyan Chen, and Zhenyu Guo.
\newblock A fractional magnetic {H}ardy-{S}obolev inequality with two
  variables.
\newblock {\em J. Math. Inequal.}, 16(1):181--187, 2022.

\bibitem{LiJiGu}
Min Liu, Fengli Jiang, and Zhenyu Guo.
\newblock Fractional {H}ardy-{S}obolev inequalities with magnetic fields.
\newblock {\em Adv. Math. Phys.}, pages Art. ID 6595961, 5, 2019.

\bibitem{MaSaVe}
Alberto Maione, Ariel~Martin Salort, and Eugenio Vecchi.
\newblock Maz{'}ya-{S}haposhnikova formula in magnetic fractional
  {O}rlicz-{S}obolev spaces.
\newblock {\em arXiv preprint arXiv:2005.04662}, 2020.

\bibitem{MiSi}
Petru Mironescu and Winfried Sickel.
\newblock A {S}obolev non embedding.
\newblock {\em Atti Accad. Naz. Lincei Rend. Lincei Mat. Appl.},
  26(3):291--298, 2015.

\bibitem{MoSk}
Kaushik Mohanta and Firoj Sk.
\newblock On the best constant in fractional $p$-{P}oincar\'e inequalities on
  cylindrical domains.
\newblock {\em arXiv preprint arXiv:2103.16845}, 2021.

\bibitem{PiSqVe}
Andrea Pinamonti, Marco Squassina, and Eugenio Vecchi.
\newblock Magnetic {BV}-functions and the {B}ourgain-{B}rezis-{M}ironescu
  formula.
\newblock {\em Adv. Calc. Var.}, 12(3):225--252, 2019.

\bibitem{Po}
Augusto~C Ponce.
\newblock An estimate in the spirit of {P}oincar{\'e}'s inequality.
\newblock {\em Journal of the European Mathematical Society}, 6(1):1--15, 2004.

\bibitem{SqVo}
Marco Squassina and Bruno Volzone.
\newblock Bourgain-{B}r\'{e}zis-{M}ironescu formula for magnetic operators.
\newblock {\em C. R. Math. Acad. Sci. Paris}, 354(8):825--831, 2016.

\bibitem{Temgoua}
Remi~Yvant Temgoua.
\newblock On the {$ s $}-derivative of weak solutions of the {P}oisson problem
  involving the regional fractional {L}aplacian.
\newblock {\em Commun. Pure Appl. Anal.}, 22(1):228--248, 2023.

\bibitem{Warma16}
Mahamadi Warma.
\newblock The fractional {N}eumann and {R}obin type boundary conditions for the
  regional fractional {$p$}-{L}aplacian.
\newblock {\em NoDEA Nonlinear Differential Equations Appl.}, 23(1):Art. 1, 46,
  2016.

\bibitem{ZuLo}
Jiabin Zuo and Juliana~Honda Lopes.
\newblock The {K}irchhoff-type diffusion problem driven by a magnetic
  fractional {L}aplace operator.
\newblock {\em J. Math. Phys.}, 63(6):Paper No. 061505, 14, 2022.

\end{thebibliography}
	\bibliographystyle{plain}
	
\end{document}